\documentclass[11pt,psfig]{article}
\usepackage{amssymb,amsmath,amscd}
\usepackage[all]{xy}
\usepackage{mathtools}
\newtheorem{lem}{Lemma}[section]
\newtheorem{thm}{Theorem}[section]
\newtheorem{cor}{Corollary}[section]
\newtheorem{prop}{Proposition}[section]

\newcommand{\commentout}[1]{}

\newcommand{\mc}{\mathcal}
\newcommand{\arr}[1]{\left( \begin{array}{clcr} #1 \end{array} \right)}

\newcommand{\diag}{{\,\rm diag}}

\newcommand{\sgn}{\, {\rm sgn}}

\newcommand{\vertiii}[1]{{\left\vert\kern-0.25ex\left\vert\kern-0.25ex\left\vert #1 
    \right\vert\kern-0.25ex\right\vert\kern-0.25ex\right\vert}}

\begin{document}
\title{ On the Growth of Rankin-Selberg $L$-functions for  $SL(2, \mathbb R)$}
\author{Hongyu He  \\
Department of Mathematics \\
Louisiana State University \\
email: hongyu@math.lsu.edu\\
}
\date{}
\maketitle

\abstract{In this paper, we establish bounds of the Rankin-Selberg $L$-function for $SL(2)$ using the supnorm of the Eisenstein series and a purely representation theoretic index over the real group. Consequently, we obtain a subconvexity bound $L(\frac{1}{2}+ it, f_1 \times f_2) \leq C (1+ |t|)^{\frac{5}{6}+\epsilon}$ for two Maass cusp forms of $SL(2, \mathbb Z)$.}

\section{Introduction}
Determining the growth rate of $L$-functions in the critical strip is an important problem in the theory of $L$-functions. Perhaps the most well-known example of such problem is the Lindelof hypothesis, that is, the Riemann-zeta function
$$|\zeta(\frac{1}{2}+ it )| \leq C_{\epsilon} (1+ |t|)^{\epsilon}$$
for every $\epsilon >0$. In the modern theory of automorphic $L$-function, generalized Lindelof hypothesis  says that all automorphic $L(\frac{1}{2}+ it)$ will grow slower than $(1+ |t|)^{\epsilon}$ (\cite{sarnak} \cite{is}). At this moment, the full Lindelof hypothesis is out of reach. A great amount of efforts has been made to establish the subconvex bounds. See \cite{fried} \cite{sarnak} \cite{mun} and the reference therein. Recall that by the Phragmen-Lindelof principle, the functional equation  implies that, a degree $n$ $L$-function has the convex bound
$$|L(\frac{1}{2}+ it)| \leq C_{\epsilon} (1+| t|)^{\frac{n}{4}+\epsilon}.$$
Any bounds that grow slower that the convex bound is called a subconvex bound. The way we stated the growth problem is often referred as the  $t$-aspect. \\
\\
The Rankin-Selberg $L$-function of $GL(n)$ is a degree $n^2$ $L$-function attached to the product of two (Adelic) cuspidal automorphic representations of $GL(n)$, often denoted by   $L(s, \pi_1 \times \pi_2)$ (\cite{cogdell}). This $L$-function can be expressed as a quotient of two Rankin-Selberg integrals, one on the automorphic side and one over $\mathbb R$. Our main result in this paper will provide evidence that $L(s, \pi_1 \times, \pi_2)$ can be  controlled by two factors, one involving the growth of the mirabolic Eisenstein series in the $t$ aspect, the other a purely representation theoretic invariant over the real group.\\
\\
In this paper, we shall treat degree 4 Rankin-Selberg $L$-functions, over the group $SL(2, \mathbb R)$. The $L$-function treated here does not need to have the usual functional equation, nor Euler product (\cite{sarnak}). Our setting and scope will be similar to those of Borel (\cite{borel}) and Harish-Chandra (\cite{hc}). 
To start, 
let $\Gamma$ be a non-uniform lattice in $SL(2)$. By a cuspidal automorphic representation of $SL(2)$, we mean an irreducible  unitary sub representation of $L^2(SL(2)/\Gamma)$ such that its constant term vanishes for all cusps (\cite{hc} \cite{la}). We say that an automorphic representation is of type $\pi$ if the automorphic representation is infinitesimally equivalent to the irreducible  smooth representation $\pi^{\infty}$.
Cuspidal representations are either in the unitary principal series, complementary series or limits of discrete series (\cite{knapp}). Unless otherwise stated, all $SL(n)$, $GL(n)$ in this paper will be over the real numbers.\\
\\
Rankin-Selberg $L$-function originated in the papers by Rankin and Selberg (\cite{rankin} \cite{selberg}). One of its most important features is that it possesses an integral representation, roughly, it can be represented by a quotient
$$L(s, \pi_1 \times \pi_2)=\frac{Tr_{aut}^{RS}(\pi_1, \pi_2, E(s, *))}{Tr_{\mathbb R}^{RS}(\pi_1, \pi_2, \mathbf 1_{s})}.$$
Here $\mathbf 1_{s}$ is the constant function $1$ on the maximal compact group $K=SO(2)$, taken as the spherical vector in the principal series representation $\mc P(2s-1,+)$ and $E(s, *)$ is the Eisenstein series.
\commentout{
Here $Tr_{aut}^{RS}$ is an integration of the product of a $\pi_1$ automorphic form, the complex conjugate of a $\pi_2$ automorphic form and the (minimal) Eisenstein series. It is a trilinear form at the representation level. The Archimedean intergal $Tr_{\mathbb R}^{RS}(\pi_1, \pi_2, \mathbf 1_{K})$ will be an integral of Whitakker functions of $\pi_1$ and $\pi_2$ against the spherical vector in the principal series $\mc P(2s-1)$ (\cite{jacquet}). }
In the literature, the integral representation of $L$-functions is  defined by choosing special testing functions, often the lowest $K$-types in $\pi_1$ and $\pi_2$. Then one can compute the Rankin-Selberg integrals and arrive at a formula of $L(s, \pi_1 \times \pi_2)$ fast. The functional equation immediately follows from the functional equation of the Eisenstein series. See for example \cite{go} \cite{cogdell}  and the references therein. However, one also lost some sights about the $L$-function in this process. The purpose of this paper is to study the Rankin-Selberg integrals at the unitary representation level and gain insights about the size of $L$-function.\\
\\
Fix the usual Iwasawa decomposition $SL(2)=KAN$. Let us state our main results.
  The first result we proved in the paper is a direct consequence of our earlier investigation in \cite{hesl2}.
\begin{thm}[A]
Let $\Gamma$ be a nonuniform lattice of $SL(2)$ such that $\Gamma \cap N=N_{\mathbb Z}$.
Let $\mc H_{\pi_1}$ and $\mc H_{\pi_2}$ be two cuspidal automorphic representations in $L^2(G/\Gamma)$.
Let $\Re(s) \in (0,1)$. The Rankin-Selberg integral 
$$Tr_{aut}^{RS}(\Phi_{f_1}, \Phi_{f_2}, E(s, \phi, g))=\int_{G/\Gamma} \Phi_{f_1}(g) \overline{\Phi_{f_2}(g)} E(s, \phi, g) d g$$
converges absolutely for $(\Phi_{f_1}, \Phi_{f_2}, \phi) \in \mc H_{\pi_1} \times \mc H_{\pi_2} \times C^{\infty}(K)$.
It yields a $SL(2)$-equivariant trilinear form on
$\mc H_{\pi_1} \times \mc H_{\pi_2} \times E(s, *, g)^{\infty}$.
This trilinear form is continuous on the first and 2nd entries, with respect to their Hilbert norms.
\end{thm}
Here $E(s, \phi,g)$ is the Eisenstein series constructed from the smooth vectors in the principal series $\mc P(2s-1,\pm)$ ( \cite{hc} \cite{la} \cite{lang} \cite{borel} \cite{knapp}). Details of this is given in Section 4. \\
\\
In order to properly define Rankin-Selberg $L$-function beyond the discrete series representations, we will need to work with $GL(2)$ to utilize the multiplicity one theorems. However, $GL(2)$ is not a semisimple group and its $L^2$ theory is more complicated due to its noncompact center. We choose to work with 
$$SL^{\pm}(2)=\{ \det g= \pm 1 \mid g \in GL(2) \}.$$
A unitary principal series or complementary series representation $\pi$ of $SL(2)$, becomes an irreducible representation of $SL^{\pm}(2)$ by defining an auxiliary action of $I_{-}=\diag(-1,1)$ (See Section 3). There are essentially two ways to do this. We call them $\pi_e$ and $\pi_o$. In the same fashion, automorphic forms of $SL(2)$ for these types, naturally extend to automorphic forms of $SL^{\pm}(2)$ in two ways which differ by a factor $\det g$. If $\pi$ is a discrete series representation or its limit, multiplicity one theorem for $SL(2)$ holds. Rankin-Selberg $L$-function can be defined over $SL(2)$ as it appeared in \cite{rankin} \cite{selberg}.\\
\\
Let us now define the Archimedean  Rankin-Selberg trilinear form over $SL^{\pm}(2)$ (\cite{jacquet}). The Iwasawa decomposition of $SL^{\pm}(2)$ is denoted by $kan$ with $k \in SO(2)$ and $a \in \mathbb R^{+} \cup \mathbb R^{-}$. For any $f_i \in \pi_i^{\infty}$ and $\phi \in C^{\infty}(K)$, define
$\phi_s(kan)=\phi(k) |a|^{-2s} $ and 
$$Tr_{\mathbb R}^{RS}(f_1, f_2, \phi_s)=\int_{G/N} Wh_{f_1}(ka) \overline{Wh_{f_2}(ka)} \phi(k) |a|^{-2s}  |a|^2 \frac{ d a}{a} d k.$$
Here $\pi_i^{\infty}$ are the smooth vectors of the irreducible unitary representation $\pi_i$ and $Wh$ stands for the Whitakker model. The Archimedean Rankin-Selberg integral only converges absolutely for $s$ in a right half plane. It can be extended by meromorphic continuation to the whole complex plane. We call this Archimedean Rankin-Selberg trilinear form. We shall remark that Miller and Schmid also studied invariant trilinear forms for $SL(2)$, in the smooth category (\cite{ms} \cite{ms0}). Their trilinear form
 may differ from $Tr_{\mathbb R}^{RS}$ perhaps by a quotient of $\Gamma$-factors.\\
\\
Define $$\mc I_{\pi_1, \pi_2}(s)=\sup \{|Tr_{\mathbb R}^{RS}(f_1, f_2, {\mathbf 1}_s) \,\,  : \,\,\, \|f_1\|=\|f_2\|=1 \}.$$
In this crude fashion,  $\mc I_{\pi_1, \pi_2}(s)$ may not be finite.  However, if $\pi_1$ and $\pi_2$ come from cuspidal representations, by multiplicity one property (\cite{loke}), $Tr_{\mathbb R}^{RS}$ differs from $Tr_{aut}^{RS}$ by a constant multiple. By Theorem $A$,
with $s$ in the critical strip, i.e. $\Re(s) \in (0,1)$,   the  induced bilinear form $TR^{RS}_{\mathbb R}(*,*, {\mathbf 1})$ on $\pi_1 \times \pi_2$ is continuous, thus bounded. Hence $\mc I_{\pi_1, \pi_2}(s)$ is finite. In order to compute $\mc I_{\pi_1, \pi_2}(s)$, we will give an intrinsic proof that $Tr_{\mathbb R}^{RS}$ is well-defined and continuous for $ \pi_1 \times \pi_2$ when $\pi_1$ and $\pi_2$ are spherical unitary principal series  and $s$ in the critical strip. We expect 
this to be the case for all nontrivial irreducible unitary representations of $SL^{\pm}(2)$.\\
\\
Our second result can be stated as follows.
\begin{thm}[B]
Let $\Gamma$ be a nonuniform lattice of $SL(2)$ such that $\Gamma \cap N=N_{\mathbb Z}$ and $G=SL^{\pm}(2)$.
Let $\mc H_{\pi_1}$ and $\mc H_{\pi_2}$ be two irreducible cuspidal representations of $G/\Gamma$. Then the 
Rankin-Selberg $L$ function  satisfies
$$|L(s, \mc H_{\pi_1} \times \mc H_{\pi_2}) | \leq C_{v, \mc H_{\pi_1}, \mc H_{\pi_2}} \frac{\sup_G(|E(s, g) v(g)|)}{\mc I_{\pi_1, \pi_2}(s)} \qquad (\forall \,\,  \Re(s) \in (0,1))$$
where $v(g)$ is a continuous bounded positive function on the fundamental domain $\mc F$ such that $v(ka_i n_i)=|a_i|^{2-\epsilon}$ for a fixed  $\epsilon \in (0, \min(2\Re(s), 2- 2\Re(s))$ on every Siegel set $S_i$.
\end{thm}
Here we decompose the fundamental domain $\mc F$ into a disjoint union of Siegel sets $S_i$ and a compact set $\mc K$ as in \cite{borel}. $G=kA_i^{\pm} N_i$ is the Iwasawa decomposition with respect to the Siegel set $S_i$. A stronger result is given in Theorem \ref{bound1}.  Theorem $A$ and $B$ are expected to hold for $SL^{\pm}(n)$ (\cite{hegl}).\\
\\
For congruence subgroups, the supnorms on the Eisenstein series are available due to recent works of  Young, Huang-Xu, Assing and Nordentoft (\cite{young}, \cite{hx} \cite{as} \cite{no}). We are left with the task of determining $\mc I_{\pi_1, \pi_2}(s)$, a purely representation theoretic invariant. This index is essentially the norm of the induced bilinear form in a proper sense.  In this paper,
we give some detailed calculation of $Tr_{\mathbb R}^{RS}$ for spherical unitary principal series $(\pi_1, \pi_2)=(\mc P(i \lambda_1,+)_e, \mc P(i \lambda_2, +)_e)$.  The precise value of $\mc I_{\pi_1, \pi_2}(s)$ is expressed in Theorem \ref{isup} as the sup norm of the Fourier coefficients of a function on the unit circle, but remains elusive numerically. Nevertheless, we prove a lower bound 
 $$\mc I_{\pi_1, \pi_2}(s) \geq  C |s|^{-1+\Re(s)} \qquad (\Re(s) \in (\frac{1}{2}, 1)).$$
 It is a little intriguing that our method has a natural barrier on the critical line $\Re(s)=\frac{1}{2}$.
The third result we prove is about the growth rate of $L(s, \mc H_{\pi_1} \times \mc H_{\pi_2})$ for congruence subgroups.
\begin{thm}[C]
Let $\Gamma$ be a congruence subgroup of $SL(2, \mathbb Z)$ or a conjugate of such a group with $\Gamma \cap N= N_{\mathbb Z}$. Let 
$(\pi_1, \pi_2)=
(\mc P(i \lambda_1, +)_{e}, \mc P(i \lambda_2, +)_{e})$. Let $\mc H_{\pi_1}$ and 
$\mc H_{\pi_2}$ be two cuspidal automorphic representation of $SL^{\pm}(2)/\Gamma$. Then for $s_0 > \frac{1}{2}$ and $\epsilon>0$, we have
$$L(s_0+i t, \mc H_{\pi_1} \times \mc H_{\pi_2}) \leq  C_{ s_0, \epsilon, \mc H_{\pi_1}, \mc H_{\pi_2}} 
(1+\|t\|)^{\frac{11}{8}-s_0+\epsilon}  $$
If $\Gamma=SL(2, \mathbb Z)$, then 
$$L(s_0+i t, \mc H_{\pi_1} \times \mc H_{\pi_2}) \leq  C_{ s_0, \epsilon, \mc H_{\pi_1}, \mc H_{\pi_2}} 
(1+\|t\|)^{\frac{4}{3}-s_0+\epsilon} . $$ 
\end{thm} 
By the functional equation and Phragmen-Lindelof principle, Theorem C implies that for $\Gamma=SL(2, \mathbb Z)$,
$$L(\frac{1}{2}+it, \mc H_{\pi_1} \times \mc H_{\pi_2}) \leq  C_{\epsilon, \mc H_{\pi_1}, \mc H_{\pi_2}} 
(1+\|t\|)^{\frac{5}{6}+\epsilon}.$$
This breaks the convexity bound. \\
\\
Throughout our paper, we use $c$ or $C$ as  symbolic constants and $c_{\epsilon, u}$ to indicate the dependence on $\epsilon$ and $u$. We retain $i$ as the square root of $-1$ even though $i$ will still be used as an integer index. This is common and the reader can easily tell the square root $i$ from the index $i$ within the context. The norm of complex numbers will be denoted by $|*|$ and we reserve $\| *\|$ for the norm of Hilbert spaces. Any integral without explicit bounds will be over the whole domain where the integration variable is defined. \\
\\
Most of the well-known results and terminology we used come from the standard books \cite{hc} \cite{la} \cite{jl} \cite{borel} \cite{lang} \cite{bump} \cite{jacquet} \cite{go}. We refer the readers to the these references for the explanation of notations, terminology, history, and citations. The literature on the subject of automorphic form on $SL(2)$ and its related $L$-functions is very rich. We shall refer the reader to   \cite{br} \cite{li} \cite{bhmm} \cite{mun} for more recent results and survey.  We gave citations where the precise statements and proofs were found. The reader should check on these citations for historical development.  We would welcome comments and suggestions about additional references and citations we should include. In any case, we believe our approach is new and raises some new challenges for real groups. In particular, knowledge about the index $\mc I_{\pi_1, \pi_2}(s)$ in various general setting can improve our understanding of the growth of Rankin-Selberg $L$-functions.    The main reference for the real group part are \cite{knapp}, \cite{cas} \cite{jacquet} and \cite{wa}.\\
\\
I would like to thank B. Speh and D. Jiang for getting me interested in automorphic forms. I would also like to thank J. Cogdell, B. Rubin, X. Li,  R. Munshi, E. Assing and A. Nordentoft for answering my questions. In particular, R. Munshi informed us that his research group have obtained some results toward the subconvexity bound $t^{\frac{15}{16}+\epsilon}$. 

\subsection{Setup for $SL(2)$}
Let $G=SL(2, \mathbb R)$. 
Let $$N=\{ n_t=\arr{ 1 & t  \\ 0 & 1} : t \in \mathbb R \}, \qquad N_{\mathbb Z}=\{ n_t, t \in \mathbb Z \},$$
 $$K=\{ k_{\theta}=\arr{\cos \theta & -\sin \theta \\ \sin \theta & \cos \theta}: \theta \in \mathbb R/ 2 \pi \mathbb Z \},$$
 $$A=\{ \arr{a & 0 \\ 0 & a^{-1}} : a \in \mathbb R^+\}.$$
 The Haar measure on $A$ will be denoted by $\frac{ d a}{a}$. 
 Let $\overline{N}$ be the opposite nilradical.
Let $M= \{ \pm I \} \subseteq K$. Fix $P=MAN$ the minimal parabolic subgroup and $\overline{P}=MA\overline{N}$ the opposite parabolic subgroup. 
Fix $ d k= \frac{1}{2 \pi} d \theta $ the invariant probability measure on $K$. We write $g=kan=k_{\theta} a n_t$ for the Iwasawa  decomposition. Let $ d g= a^2  d t \frac{d a}{ a} d k d n$ be the invariant measure of $G$ under the $KAN$ decomposition.
\\
\\
Let $\Gamma$ be a nonuniform lattice of  of $G$ such that 
$\Gamma \cap N $ is nontrivial. This is always possible by moving the cusp to $0$ using a group translation (\cite{borel}). By rescaling $N$ (\cite{borel}),  we assume that 
$\Gamma \cap N=N_{\mathbb Z} $.
 $G/\Gamma$ has a finite volume and a finite number of cusps, $z_1, z_2, \ldots, z_l$. Fix a fundamental domain $\mc F$. Write
 $\mc F$,  as a disjoint union of Siegel sets  $S_1, S_2, \ldots S_l$ and a compact set $\mc K$ (\cite{borel}).   Over each Siegel set $S_i$, the invariant measure can be written as $d g= a_i d a_i d t_i d k$. Since $\Gamma$ action is on the right, our standard Siegel set will be near $0$, 
 not $\infty$. Roughly, our results on $g$ shall translate into results on $g^{-1}$ for $\Gamma \backslash G$.\\
 \\
 We often use $\pi$ to denote the Hilbert space of a unitary representation and $\pi^{\infty}$ to denote the  space of smooth vectors of $\pi$. When $\pi$ is admissible and no obvious Hilbert structure is imposed, $\pi^{\infty}$ will denote the admissible representation in the smooth category.

\section{Supnorm of Eisenstein series} 
We follow the definition of Eisenstein series from (\cite{hc} \cite{la} \cite{lang} \cite{borel}). Let $\phi \in C^{\infty}(K)$. Let $\phi_s(kan)= \phi(k) a^{-2s}$. Define the Eisenstein series
$$E(s, \phi, g)=\sum_{\gamma \in \Gamma/\Gamma_N} \phi_s(g \gamma).$$
$E(s, \phi, g)$ converges absolutely for $\Re(s) >1$ and has a meromorphic continuation to the whole complex plane $\mathbb C$. In addition, $E(s, \phi, g)$ is holomorphic in $\Re(s) \geq \frac{1}{2}$ except finitely many simple poles in $(\frac{1}{2}, 1]$. At such a pole, $E(s, \phi, g)$ will be a square integrable automorphic function (see Theorem 11.13 \cite{borel}). \\
\\
Eisenstein series are generally slowly increasing at the cusps even though they may be rapidly decaying near some cusps. Their growth at the cusps are controlled by their constant terms (\cite{hc} \cite{la}), in our context
$$E(s, \phi, kan)= C_{s, \phi}^{(1)}(k) a^{-2s}+C^{(2)}_{{s, \phi}}(k) a^{-2+2s} +F(s, \phi, kan)$$ near the cusp zero. Here $F(s, \phi, kan)$ will be rapidly decaying near zero.\\
\\
Let $\mathbf 1$ be the constant function $1$ on $K$. We may simply write $E(s, \mathbf 1, kan)$ as $E(s, kan)$. Recently, there have been growing interests of studying the sup bounds of $E(s, kan)$ in the $s-$ aspect (\cite{young} \cite{hx} \cite{as} \cite{no}). We have
\begin{thm}[Young, Huang-Xu, Assing, Nordentoft]\label{yhxan}  Let $\Gamma$ be a congruent subgroup of 
$SL(2, \mathbb Z)$ or a      conjugate of such group. There exists constant $C_{\Gamma}$ such that
$$E_{\Gamma}(\frac{1}{2}+ it, ka_i n_i) \leq C_{\Gamma, \epsilon} a_{i}^{-1} |t|^{\frac{3}{8}+\epsilon}$$
for each Siegel set $S_i$ and $\mc K$, uniformly on $t$.
If $\Gamma=SL(2, \mathbb Z)$, this bounds can be improved to
$$E_{\Gamma}(\frac{1}{2}+ it, kan) \leq C_{\Gamma, \epsilon} a^{-1} |t|^{\frac{1}{3}+\epsilon},$$ 
on its standard fundamental domain.
\end{thm} 
The first statement is due to Young, Huang-Xu and Assing. The 2nd  statement  is due to Nordentoft. Applying Phragmen-Lindelof principle, we have

\begin{cor}\label{yhxan1}
Let $\Gamma$ be a congruence subgroup of 
$SL(2, \mathbb Z)$ or a conjugate of such group. There exists constant $C_{\Gamma}$ such that for any $\sigma \geq \frac{1}{2}$
$$E_{\Gamma}(\sigma+ it, ka_in_i) \leq C_{\Gamma, \epsilon} a_i^{-2 \sigma} |t|^{\frac{3}{8}+\epsilon}$$
on each Siegel set $S_i$ and 
for $\Gamma=SL(2, \mathbb Z)$, 
$$E_{\Gamma}(\sigma+ it, kan) \leq C_{\Gamma, \epsilon} a^{-2 \sigma} |t|^{\frac{1}{3}+\epsilon}.$$
\end{cor}
Proof:  Fix $g=k a_i n_i$ in the Siegel set $S_i$.  Consider $a_i^{2 (\sigma+ it)} E_{\Gamma}(\sigma+ it, ka_i n_i)$. This function is bounded by $C |t|^{\frac{3}{8}+\epsilon}$ on $\sigma=\frac{1}{2}$ and by $C$ on $\sigma=\frac{3}{2}$. It can have only finite number of poles between. Our corollary follows from the Phragmen-Lindelof principle. $\Box$ \\
\\ 
It is not clear whether these results can be proved for noncongruence subgroups. The group $G$ acts on the Eisenstein series from left
$$L(g)E(s, \phi, h)=E(s, \phi, g^{-1}h)= \sum_{\gamma \in \Gamma/\Gamma_{N}} \phi_s(g^{-1} h \gamma)= \sum_{\gamma \in \Gamma/\Gamma_{N}} (L(g) \phi_s)(h \gamma).$$
From representation theoretic viewpoint,  Eisenstein series is a realization of principal series $\mc P(2s-1, \pm)$ that we shall introduce.
 \section{ Principal series representations}
Let  $\chi_{-}(x)$ be the sign character on $\mathbb R - \{0 \}$ and $\chi_{+}(x)$ be the trivial character.
 For the sake of the introducing the Whitakker model,
we construct $\mc P(u, \pm)$ on $G/\overline{P}$: $f(kman)= \chi_{\pm}(m) f(k) a^{u+1}$. Let $\mc P(u, \pm)^{\infty}$ be the space of $f$ with $f(km)=\chi_{\pm}(m) f(k)$ and $f|_K$  smooth.
Consider the noncompact picture (\cite{knapp}). We have for any $g=\arr{a & b \\ c & d}$, $f \in \mathcal P(u, \pm)^{\infty}$,
$$\pi_{u, \pm}(g) f (x)=\chi_{\pm}(a-cx) |a - c x|^{-1-u} f(\frac{d x-b}{a-cx}).$$

In particular, we have
$$\pi_{u, \pm}\arr{a & 0 \\ 0 & a^{-1}} f(x)=|a|^{-1-u} f(a^{-2} x), \qquad (a \in \mathbb R^+);$$
$$\pi_{u, \pm} \arr{1 & b \\ 0 & 1} f(x)= f(x-b);$$
$$\pi_{u, \pm} \arr{\cos \theta & - \sin \theta \\ \sin \theta & \cos \theta} f(x)=\chi_{\pm}(\cos \theta- x \sin \theta) |\cos \theta- x \sin \theta|^{-1-u} f(\frac{ x \cos \theta+ \sin \theta}{\cos \theta- x \sin \theta}).$$
There is a $G$-invariant pairing between $\mc P(u, \pm)^{\infty}$ and $\mc P(-u, \pm)^{\infty}$. This allows us to write the dual space of $\mc P(u, \pm)^{\infty}$ as $\mc P(-u, \pm)^{-\infty}$.\\
\\
We will use $f \in \mc P(u, \pm)$ to denote the noncompact picture.  The space $\mc P(u, \pm)^{\infty}$ will then be a subspace of infinitely differentiable functions on $N \cong \mathbb R$ satisfying certain conditions at infinity. We will use $\tilde f \in \mc P(u, \pm)$ to denote the compact picture. In the compact picture, $\mc P(u, \pm)^{\infty}$ can be identified with $C^{\infty}(K/M, \chi_{\pm})$. Hence $\mc P(u, +)^{\infty}$ will consist of even smooth functions on the unit circle and $\mc P(u, -)^{\infty}$ will consist of odd smooth functions on the unit circle. \\
\\
For our construction of Eisenstein series, we actually model $\mc P(2s-1, \pm)$ on $G/P$. The sole purpose of doing this, as we shall see later, is to define the Rankin-Selberg integral over $\mathbb R$ correctly. 
Principal series representations of $G$ can also be constructed using homogeneous distributions on $\mathbb R^2$. See for example \cite{cas}. The advantage of this construction is that it is independent of choices of minimal parabolic $P$. Different choices  of the minimal parabolic amount to different coordinate system for the projective space $\mathbb P \mathbb R^1 \cong G/P$. \\
\\
Unitarizable principal series  include unitary principal series $\mc P(u, \pm)$ (with $u \in i \mathbb R$) and complementary series $\mc P(u,+)$ (with $u \in (-1,0) \cup (0,1)$). All of these representations are irreducible except $\mc P(0, -)$. In addition $\mc P(u, \pm) \cong \mc P(-u, \pm)$. For unitarizable principal series,  $\mathcal P(u, \pm)$ will refer to the Hilbert space with the $G$-invariant norm.  \\
\\
Now we shall turn to the disconnected group
 $$SL^{\pm}(2)=\{ \det g= \pm 1 \mid g \in GL(2, \mathbb R) \}.$$
 We must use this group in order to apply the multiplicity one theorem to construct the $L$-function. Let
 $I_{-}=\diag(-1, 1)$. Then $SL^{\pm}(2)= I_{-} SL(2) \sqcup SL(2)$. Principal series representations of $SL^{\pm}(2)$, restricted to $SL(2)$, are principal series of $SL(2)$. We can easily realize them on $\mc P(u, \pm)$ by defining the auxiliary action of $I_{-}$:
 $$\pi_{u, \pm, e} (I_{-}) f(x)= f(-x);$$
 $$\pi_{u, \pm, o}(I_{-}) f (x)=- f(-x).$$
It can be easily verified that
$$\pi_{u, \pm}(I_{-} g I_{-})=\pi_{u, \pm, *}(I_{-}) \pi_{u, \pm}(g)\pi_{u, \pm, *}(I_{-}) \qquad (g \in SL(2)).$$
Hence $ \pi_{u, \pm, *}$ is a group representation of $SL^{\pm}(2)$.
 Even though the representation space $\mc P(u, \pm)$ remains intact in this construction, we will use $\mc P(u, \pm)_{e}$ and $\mc P(u, \pm)_{o}$ to indicate the group action of $SL^{\pm}(2)$. There is not much difference between $\pi_{u, \pm, e} $ and $\pi_{u, \pm, o}$:
 $$\pi_{u, \pm, e} \cong \pi_{u, \pm, o} \otimes \det.$$ 
 \vspace{.1in}
 
\noindent 
Fix $\Gamma \in SL(2)$. We shall use $E_{\Gamma}(s, \phi, g)_{e}$ to denote the Eisenstein series associated with $\mc P(2s-1, \pm)_{e}$ and $E_{\Gamma}(s, \phi, g)_{o}$ to denote the Eisenstein series associated with $\mc P(2s-1, \pm)_{o}$. When $\Gamma$ is fixed, we will drop the subscript $\Gamma$. The Eisenstein series $E(s, \phi, g)$ with $g \in SL^{\pm}(2)$ will be read as $E(s, \phi, g)_{e}$. We have $$E_{\Gamma}(s, \phi, g)_{e}= \det(g) E_{\Gamma}(s, \phi, g)_{o}.$$
 
 \noindent
 Let $A^{\pm}= A \sqcup I_{-} A$. Then we will still have Iwasawa decomposition
 $SL^{\pm}(2)=KA^{\pm} N$. The computation on $SL^{\pm}(2)$ will be similar to computation on $SL(2)$.

\section{Rankin-Selberg Trilinear form: Automorphic side}
Let $\Gamma$ be a nonuniform lattice in $SL(2)$ and $G=SL^{\pm}(2)$. Then $G/\Gamma$ has a finite volume and a finite number of cusps, $z_1, z_2, \ldots, z_l$. Fix a fundamental domain $\mc F$ for 
 $G/\Gamma$ and write it as the disjoint union of Siegel sets  $S_1, S_2, \ldots S_l$ with a compact set $\mc K$ (\cite{borel}). Let $ d g= |a|\,  d a \,d t \, d k$ be the invariant measure of $G$ under the $KA^{\pm} N$ decomposition. Over each Siegel set $S_i$, the invariant measure will be written as $d g=| a_i| d a_i d t_i d k$. We have
 
 \begin{thm}[\cite{hesl2}]\label{main1} Let $\Gamma$ be a nonuniform lattice in $SL(2, \mathbb R)$. Let $G= SL^{\pm}(2)$. Let $\mathcal H \subseteq L^2(G/\Gamma)$ be a cuspidal automorphic representation of type $\pi$. Let $\epsilon >0$. Given any positive measure $\nu$ on $G/\Gamma$ such that $\nu$ is bounded by $ dg$ on $\mc K$ and bounded by $|a_i|^{\epsilon} \frac{d a_i}{|a_i|} d t_i d k$ on $S_i$,  there exists a constant $C$ depending on $\nu$ (hence on $\epsilon$) and $\mc H$ such that

$$ \| f \|_{L^2(G/\Gamma, \nu)} \leq C \|f \|_{L^2(G/\Gamma, d g)}, \qquad ( f \in \mc H).$$
\end{thm}
In Section 3 \cite{hesl2}, this theorem is proved for $SL(2)$ for each Siegel set. The $SL^{\pm}(2)$ case is not much different. $\Box$ \\
\\
Let $\pi_1$ and $\pi_2$ be two irreducible unitary representations of $G=SL^{\pm}(2)$. Let $\pi_1^{\infty}$ and $\pi_2^{\infty}$ be the Frechet spaces of smooth vectors.  Suppose that $\Gamma \cap N =N_{\mathbb Z}$. Let  $\mc H_{\pi_1}$ and $\mc H_{\pi_2}$ be two cuspidal representations of type $\pi_1$ and $\pi_2$ in $L^2(G/\Gamma)$. For each $f \in \pi$, let $\Phi_{f}$ be the corresponding automorphic function in $L^2(G/\Gamma)$ under a fixed identification of $\pi$ with $\mc H_{\pi}$. The identification is unique upon a normalization of Fourier coefficient, for example, by making $a_1=1$. For  $f_i \in \pi_{i}^{\infty}$, $\phi \in C^{\infty}(K)$,   define the Rankin-Selberg trilinear form:
 $$Tr_{aut}^{RS}(\Phi_{f_1}, \Phi_{f_2}, E(s, \phi, *)_{e, o})=\int_{G/\Gamma} \Phi_{f_1}(g) \overline{\Phi_{f_2}(g)} E(s, \phi, g)_{e, o} d g.$$
   We then obtain a trilinear form, complex linear on the first and third entries, complex conjugate linear on the second entry. We may call the right hand side Rankin-Selberg integral. Since $\Phi_{f_i}$ are rapidly decaying near the cusps and the Eisenstein series grow slower than a polynomial near the cusps, Rankin-Selberg integral converges absolutely. From the analytic properties of $E_{\Gamma}$ (\cite{borel}), we have 
 \begin{prop}
 The trilinear form $Tr_{aut}^{RS}$ is well defined for all  triple  in $\mc H_{\pi_1}^{\infty} \times \mc H_{\pi_2}^{\infty} \times E(s, *, *)^{\infty}_{e,o}$ as long as the Eisenstein series $E(s, *,*)$ is well-defined. As $s$ varies, the trilinear form has a meromorphic continuation on the right half plan $\Re(s) \geq 0$. It has only a finite number of simple poles in the interval $(\frac{1}{2}, 1 ]$ in the region  $\Re(s) \geq \frac{1}{2}$. 
 \end{prop}
 
Here $E(s, *, *)^{\infty}_{e,o}$ denote the  space of Eisenstein series $E(s, \phi, g)_{e,o}$ with $\phi \in C^{\infty}(K)$. It corresponds to $\mc P(2s-1, \pm)^{\infty}_{e,o}$ depending on whether $\phi \in C^{\infty}(k/M)$ or $C^{\infty}(K/M, \chi_{-})$. By Theorem \ref{main1}, the Rankin-Selberg integral actually converges for all functions in $\mc H_{\pi_1} \times \mc H_{\pi_2} \times E(s, *, *)^{\infty}_{e,o}$ when $2\Re(s)-1 \in (-1,1)$.  

\begin{thm}\label{main2}
Let $\Re(s) \in (0,1)$. The Rankin-Selberg integral converges absolutely for $(\Phi_{f_1},\Phi_{f_2}, (\phi_s)_{e,o}) \in \mc H_{\pi_1} \times \mc H_{\pi_2} \times \mc P(2s-1,\pm)_{e, o}^{\infty}$ as long as the Eisenstein series $E(s, \phi, *)$ is well-defined.
It yields a $G$-equivariant trilinear form on
$\mc H_{\pi_1} \times \mc H_{\pi_2} \times \mc P(2s-1,\pm)_{e, o}^{\infty} $.
This trilinear form is continuous on the first and 2nd entries.

\end{thm}
 Continuity of $Tr_{aut}^{RS}$ on the third entry can be proved, but is not needed in this paper. So we choose not to include it in our theorem. We remark that this extension is valid exactly when $s$ is in the critical strip.\\
\\
Proof: Fix $s$ with $\Re(s) \in (0,1)$. Let $(\Phi_{f_1},\Phi_{f_2}, (\phi_s)_{e,o}) \in \mc H_{\pi_1} \times \mc H_{\pi_2} \times \mc P(2s-1, \pm)_{e,o}$. We consider the measure $\nu= |E(s, \phi, g)| d g$ on a fundamental domain $\mc F$. Over each Siegel set $S_i$, $|E(s, \phi, g)|$ is bounded by a constant multiple of $|a_i|^{-2\Re(s)}+|a_i|^{-2+2\Re(s)}$. Since $\Re(s) \in (0,1)$ and 
$d g= |a_i|^2 \frac{ d a_i}{|a_i|} d k d n_i$ ,  $\nu$ will be bounded by a constant multiple of $|a_i|^{\epsilon} \frac{ d a_i}{|a_i|} d k d n_i$ on $S_i$ for some $\epsilon >0$. In addition, over the compact set  
$\mc K$, $|E(s, \phi, g)|$ is a bounded and smooth. $\nu$ will be bounded by a multiple of $d g$ on $\mc K$.  We then apply Theorem \ref{main1}. We have $ \Phi_{f_1}, \Phi_{f_2} \in L^2(\mc F, \nu)$. By the Cauchy-Schwartz inequality, 
$Tr_{aut}^{RS}(\Phi_{f_1}, \Phi_{f_2}, E(s, \phi, g)_{e,o})$ converges absolutely. In addition, it is continuous  on $\Phi_{f_1}$ and $\Phi_{f_2}$. $\Box$

\section{Rankin-Selberg Trilinear form over $\mathbb R$}
We shall now define the local analogue of Rankin-Selberg trilinear form over $\mathbb R$ (\cite{jacquet}). 
 Fix two irreducible infinite dimensional unitary representation $\pi_1$ and $\pi_2$ of $G=SL^{\pm}(2)$. Let $i=1,2$. We  use $\pi_i$  to also denote the Hilbert space and reserve $\mc H_{\pi_i}$ for the Hilbert space of the cuspidal representation. Let $\pi_i^*$ be the contragradient representation. We use $\langle *, * \rangle$ to denote the canonical $G$-equivariant nondegenerate pairing. The space of continuous linear functionals on $\pi_i^{\infty}$, will be denoted by
$(\pi_i^*)^{-\infty}$.\\
\\
 Let $\psi_i$ be the unique Whitakker vector in $(\pi_i^*)^{-\infty}$. Later we shall specify the Whittaker vector when we restrict our attention to principal series representations.\\
\\
For any $f_i \in \pi_i$, define the Whitakker function
$$Wh_{f_i}(g)=\langle f, \pi^*(g) \psi_i \rangle.$$
It is not obvious that $Wh_{f_i}(g)$ is a locally integrable function on $G$. Nevertheless, if $f_i \in \pi_i^{\infty}$, $Wh_f(g)$ will be a smooth function on $G$ satisfying
$$Wh_{f_i}(gn_t)=Wh_{f_i}(g) \exp 2 \pi i t.$$
To make sure we have the right kind of covariance, we verify that
$$L(g) Wh_{f}(h)= Wh_{f}(g^{-1} h)=\langle f, \pi^*(g^{-1} h) \psi \rangle= \langle \pi(g) f, \pi^*(h) \psi \rangle=Wh_{\pi(g) f}(h).$$
Let $\chi_{e}(a)=1$ and $\chi_{o}(a)=-1$ if $ a <0$; $\chi_{e,o}(a)=1$ if
$a>0$ .
For any $f_i \in \pi_i^{\infty}$ and $\phi \in C^{\infty}(K/M, \chi_{\pm})$, define
$\phi_s(kan)_{e,o}=\phi(k) |a|^{-2s} \chi_{e,o}(a)$ and 
$$Tr_{\mathbb R}^{RS}(f_1, f_2, (\phi_s)_{e,o})=\int_{G/N} Wh_{f_1}(ka) \overline{Wh_{f_2}(ka)} \phi(k) |a|^{-2s} \chi_{e,o}(a) |a|^2 \frac{ d a}{a} d k.$$
Observe that   $Wh_{f_1}(ka) \overline{Wh_{f_2}(ka)}$ is a smooth function on $G/N$ and $|a|^{-2s} \chi_{e,o}(a) \phi(k)$ can be realized as a smooth vector in the principal series $\mc P(2s-1,\pm)_{e,o}$ modeled on $G/P$. We see that $Tr_{\mathbb R}^{RS}$ is a $G$-equivariant trilinear form on
$$\pi_1^{\infty} \times \pi_2^{\infty} \times \mc P(2s-1, \pm)^{\infty}_{e,o}.$$
To simplify our notation, $\phi_s$ will denote ${\phi_s}(kan)_{e}$. 
Now we need to know the convergence of $Tr_{\mathbb R}^{RS}$. The following is well-known.

\begin{prop} Let $\pi_1$ and $\pi_2$ be two irreducible unitary representations of $G$. Suppose that
$f_i \in \pi_i^{\infty}$ and $\phi \in C^{\infty}(K)$. Then  $Tr_{\mathbb R}^{RS}(f_1, f_2, \phi_s)$ converges absolutely for all $\Re(s) >1$. If $\pi_1$ and $\pi_2$ are tempered, then $Tr_{\mathbb R}^{RS}(f_1, f_2, \phi_s)$ converges absolutely for $\Re(s) >0$.
\end{prop}

Later we shall examine $Tr_{\mathbb R}^{RS}$ in more details. In the literature, when this trilinear form is applied, it often include its meromorphic continuation to $s \in \mathbb C$. This makes sense only when $f_1, f_2, \phi_s$ are smooth. The convergence becomes less of an issue in the smooth category.  Similar trilinear forms were studied in \cite{loke} \cite{ms}. \\
\\
For the purpose of this paper, we must examine  when $Tr_{\mathbb R}^{RS}$ is well-defined for the full unitary representations $\pi_1$ and $\pi_2$. This will be done in Section 7. It is also desirable to find the maximal value of $Tr^{RS}_{\mathbb R}
(f_1, f_2, \phi_s)$ for $\|f_1\|=\|f_2\|=1$ and $\phi_s$ fixed, which turns out to be hard. We will focus our attention on $\phi_s= \mathbf 1_s$ in Section 8. The corresponding integral on the automorphic side will have the (spherical) Eisenstein series $E(s, g)$ as the third entry. This is precisely what Rankin and Selberg used to define their $L$-functions. Technically, choosing the spherical vector $\phi=\mathbf 1$ simplifies $Tr_{\mathbb R}^{RS}$ considerably.  

\section{Integral Representation of Rankin-Selberg $L$-function}
One of the most important features of Rankin-Selberg $L$-function is its integral representation.  The adelic language presentation of such theory for $GL(n)$ can be found in \cite{cogdell}. For $G=SL^{\pm}(2, \mathbb R)$, Rankin-Selberg $L$-function can be  defined for a much larger class of discrete subgroups.  we shall give a brief review of Rankin-Selberg $L$-function for $G=SL^{\pm}(2)$. The main references will be \cite{bump} and \cite{go}. Unless otherwise stated, Eisenstein series $E(s, \phi, g)$ will refer to $E(s, \phi, g)_{e}$.\\
\\
Let $\Gamma$ be a nonuniform lattice in $SL(2)$ with 
$ \Gamma \cap N=N_{ \mathbb Z }$.
For two cuspidal automorphic representations $\mc H_{\pi_1}$ and $\mc H_{\pi_2}$ of $G/\Gamma$
 we define formally.
 \begin{equation}\label{l(s)}
 L(s, \mc H_{\pi_1} \times \mc H_{\pi_2})=\frac{Tr_{aut}^{RS}(\Phi_{f_1}, \Phi_{f_2}, E(s, \phi, *))}{Tr_{\mathbb R}^{RS}(f_1, f_2, \phi_s)}
 \end{equation}
 for any $\Phi_{f_1} \in \mc H_{\pi_1}^{\infty}, \Phi_{f_1} \in \mc H_{\pi_1}^{\infty}, \phi_s \in \mc P(2s-1,+)^{\infty}_{e}$ as long as the denominator is not zero. \\
 \\
 This definition is valid for all cuspidal representations of $G/\Gamma$, due to the multiplicity one theorem for $G$, namely, the $G$-equivariant  trilinear form 
 $$Tr: \pi_1^{\infty} \times \pi_2^{\infty} \times \mc P(2s-1, +)^{\infty}_{e, o} \rightarrow \mathbb C$$
 is unique up to a constant (\cite{loke}). Since the multiplicity one theorem fails for $SL(2)$,  we introduce the group $G=SL^{\pm}(2)$ at the beginning. There are other remedies besides introducing $SL^{\pm}(2)$. One remedy is to follow \cite{bump} and \cite{go} and to define even and odd automorphic forms. Another remedy is to go to the group $GL(2)$ and introduce a center and an extra parameter for the representation. These two approaches will likely make our results  more complicated to state and to prove. \\
 \\
Our definition of $L$-function depends on $\Phi$: the identification of $\pi_i$ with $\mc H_{\pi_i}$. This can be easily resolved by normalizing the Fourier-Whitakker coefficients $\{a_n \}$ and $\{b_n \}$ of the automorphic forms properly.\\
\\
 The validation of our definition, that is, the quotient is indeed an $L$-function, follows from the standard unwinding process of the Rankin-Selberg integral $Tr_{aut}^{RS}$. For Maass cusp forms, $L(s, \mc H_{\pi_1} \times \mc H_{\pi_2})$ will be of the form $\sum_{n \neq 0} \frac{a_n \overline{ b_n}}{|n|^s}$ or of the form $\sum_{n \neq 0} \sgn(n) \frac{a_n \overline{ b_n}}{|n|^s}$.  This is true at least for any $\Re(s)>1$ by the unwinding process.  \\
\\
Secondly,  the trilinear form  $Tr_{aut}^{RS}$ is well-defined for $\Re(s) >1$ and is analytic on the variable $s$ except where the Eisenstein series $E_{\Gamma}(s, *,*)$ has poles.  $Tr_{\mathbb R}^{RS}$ is also well defined for $\Re(s)>1$ and has analytic continuation over $\mathbb C$ except some poles.
For a difference choice $(f_1^{\prime}, f_2^{\prime}, \phi^{\prime})$,
$$\frac{Tr_{aut}^{RS}(\Phi_{f_1}, \Phi_{f_2}, E(s, \phi, *))}{Tr_{\mathbb R}^{RS}(f_1, f_2, \phi_s)}=\frac{Tr_{aut}^{RS}(\Phi_{f_1^{\prime}}, \Phi_{f_2^{\prime}}, E(s, \phi^{\prime}, *))}{Tr_{\mathbb R}^{RS}(f_1^{\prime}, f_2^{\prime}, \phi^{\prime}_s)}$$
for $\Re(s)>1$, hence for all $s \in \mathbb C$. We see that
$L(s, \mc H_{\pi_1} \times \mc H_{\pi_2})$ does not depend on the choices of $(f_1, f_2, \phi)$. In addition, the quotient $L(s, H_{\pi_1} \times \mc H_{\pi_2})$ will be analytic on $\Re(s) \geq \frac{1}{2}$ except possibly a finite number of poles from Eisenstein series $E_{\Gamma}(s, \phi, g)$ (\cite{borel}).
\begin{thm}\label{bound1} Rankin-Selberg $L$ function with $\Re(s) \in (0,1)$
$$|L(s, \mc H_{\pi_1} \times \mc H_{\pi_2}) | \leq C_{v, \mc H_{\pi_1}, \mc H_{\pi_2}} \inf\{\frac{\sup_G(|E(s, \phi, g) v(g)|)}{|Tr_{\mathbb R}^{RS}(f_1, f_2, \phi_s)|} \,\,  : \,\,\, \|f_1\|=\|f_2\|=1, \phi \in C^{\infty}(K/M) \}$$
where $v(g)$ is a bounded continuous positive function on $\mc F$ such that $v(ka_i n_i)=|a_i|^{2-\epsilon}$ for a fixed  $\epsilon \in (0, \min(2\Re(s), 2- 2\Re(s))$ on every Siegel set $S_i$.
\end{thm}
Here the constant $C_{v, \mc H_{\pi_1}, \mc H_{\pi_2}}$ can be determined by using certain $L^2$-norm of Fourier coefficients $\{a_i \},\{b_i \}$ as in \cite{hesl2} upon a proper choice of $v(g)$. A more complete and complicated formula can be obtained by taking the inf of all possible choices of $\epsilon$. We may not gain much from this because the constant $C_{v, \mc H_{\pi_1}, \mc H_{\pi_2}}$ may also change. However if the choice of $\epsilon$ and in turn $v$ will alter $\sup_G(|E(s, \phi_s, g) v(g)|)$ in the $s-$aspect, we will gain improvement of bounds on $L(s, \mc H_{\pi_1} \times \mc H_{\pi_2})$.\\
\\
Proof: Let $\nu=v(g)^{-1} d g$. Then $\nu=|a_i|^{\epsilon} \frac{d a_i}{a_i} d n_i d k$ on each Siegel set $S_i$. It is bounded above by 
$\frac{1}{\min( v(g), g \in \mc K)} dg $ on $\mc K$. By Cauchy-Schwartz inequality and Theorem \ref{main1}, we have
\begin{equation}
\begin{split}
|L(s, \mc H_{\pi_1} \times \mc H_{\pi_2})|= & \frac{|Tr_{aut}^{RS}(\Phi_{f_1}, \Phi_{f_2}, E(s, \phi, *))|}{|Tr_{\mathbb R}^{RS}(f_1, f_2, \phi_s)|} \\
\leq & \frac{ \| \Phi_{f_1}(g) v(g)^{-\frac{1}{2}} \|_{L^2(G/\Gamma)} \| \Phi_{f_2}(g) v(g)^{-\frac{1}{2}} \|_{L^2(G/\Gamma)} \|E(s, \phi, g) v(g) \|_{\sup} }{|Tr_{\mathbb R}^{RS}(f_1, f_2, \phi_s)|} \\
= & \frac{ \| \Phi_{f_1}(g)  \|_{L^2(G/\Gamma, \nu)} \| \Phi_{f_2}(g)  \|_{L^2(G/\Gamma, \nu)} \|E(s, \phi, g) v(g) \|_{\sup} }{|Tr_{\mathbb R}^{RS}(f_1, f_2, \phi_s)|} \\
\leq & C_{v, \mc H_{\pi_1}, \mc H_{\pi_2}} \frac{ \| f_1 \|  \| f_2 \| \|E(s, \phi, g) v(g) \|_{\sup} }{|Tr_{\mathbb R}^{RS}(f_1, f_2, \phi_s)|}. \\
\end{split}
\end{equation}
Our theorem then follows by taking the $inf$ over $(f_1,f_2, \phi)$. $\Box$ \\
\\
The bound in Theorem \ref{bound1} may be too difficult to compute. We may fix a testing function and let $\phi=\mathbf 1$, the spherical vector. We then have
\begin{thm}\label{main3}
Rankin-Selberg $L$ function with $\Re(s) \in (0,1)$
$$|L(s, \mc H_{\pi_1} \times \mc H_{\pi_2}) | \leq C_{v, \mc H_{\pi_1}, \mc H_{\pi_2}} \frac{\sup_G(|E(s, g) v(g)|)}{\sup\{|Tr_{\mathbb R}^{RS}(f_1, f_2, {\mathbf 1}_s)|\,\,  : \,\,\, \|f_1\|=\|f_2\|=1 \}}$$
where $v(g)$ is a bounded continuous positive function on $\mc F$ such that $v(ka_i n_i)=|a_i|^{2-\epsilon}$ for a fixed  $\epsilon \in (0, \min(2\Re(s), 2- 2\Re(s))$ on  every Siegel set $S_i$.
\end{thm}
Let us define $$\mc I_{\pi_1, \pi_2}(s)=\sup \{|Tr_{\mathbb R}^{RS}(f_1, f_2, {\mathbf 1}_s) \,\,  : \,\,\, \|f_1\|=\|f_2\|=1 \}.$$ Then
for $\Re(s) \in (0,1)$, $$|L(s, \mc H_{\pi_1} \times \mc H_{\pi_2}) | \leq C_{v, \mc H_{\pi_1}, \mc H_{\pi_2}}  \frac{\sup_G(|E(s, g) v(g)|)}{\mc I_{\pi_1, \pi_2}(s)}.$$
The invariant $\mc I_{\pi_1, \pi_2}(s)$ is an invariant over the real group and $\sup$-norm is on the Eisenstein series.  $C_{v, \mc H_{\pi_1}, \mc H_{\pi_2}}$ depends only on $v$ and the $L^2$-norm of Fourier coefficients of $\mc H_{\pi_1}$ and $\mc H_{\pi_2}$. We have essentially bounded the $L$-function by three separate factors. It perhaps is not difficult to see most of our results will carry over with minimal changes to $SL^{\pm}(n)$ (\cite{hegl}. In what follows, we will will apply Theorem \ref{main3} to establish subconvex bounds for $L(s, \mc H_{\pi_1} \times \mc H_{\pi_2} )$ for congruence subgroups of $SL(2)$.\\
\\
\section{$Tr_{\mathbb R}^{RS}$ : Spherical principal series case}
Let $G=SL^{\pm}(2)$. Recall that $G=KA^{\pm} N$. We still write $g=kan$ with the understanding that $a$ can be a positive or negative number.  Since the computation of even and odd cases will be similar, we may write $\chi(a)$ with the understanding that $\chi$ could be either $\chi_{e}$ or $\chi_{o}$.  We may regard  $\pi = \mc P(i \lambda, \pm )_{e,o}$ as the  principal series $\mc P(i \lambda, \pm)$ of $SL(2)$ with an auxiliary acting of $I_{-}$. The action of $I_{-n}$ does not change the Fourier-Whitakker coefficients, Whitakker functional, or the Hilbert structure, or the action of the Lie algebra. We can still carry out our discussion  within the frame of $SL(2)$.
\\
\\
The case of spherical principal series will be  most illuminating. 
We  start with a pair of spherical unitary representations 
$(\pi_1, \pi_2)$ of $G$, where $\pi_i = \mc P(i \lambda_i, +)_{e,o}$ with an unspecified $\chi_i$. We also use $\chi_0$ to denote the action of $I_{-}$ on $(\phi_s)_{e,o}$. We shall discuss in more details the trilinear form
$$Tr_{\mathbb R}^{RS}(f_1, f_2, (\phi_s)_{e,o})=\int_{G/N} Wh_{f_1}(ka) \overline{Wh_{f_2}(ka)} \phi(k) |a|^{-2s} a^2 \chi_0(a) \frac{ d a}{a} d k$$
 Our goal is to make the trilinear form explicit enough so we can determine the range of convergence for $f_1$,$ f_2$ and 
$s$. Based on that,  we will try to give some estimate of $\mc I_{\pi_1, \pi_2}(s)$ in the next section.  \\
\\
We start with $f_i \in  \mc P(i \lambda_i, +)^{\infty}$. At this moment $Tr_{\mathbb R}^{RS}$ is only defined for $\Re(s)>1$. We shall end with a formula in terms of the compact picture $\tilde f_i$ which will demonstrate the analytic continuation of $Tr_{\mathbb R}^{RS}$ to at least $\Re(s) >0$. The Fourier transform will be denoted by $\mc F$:
$$\mc F f(\xi)= \int f(x) \exp (- 2 \pi i x \xi ) d x.$$
\subsection{Expression of $Wh_{f}(kan)$ with $a \in A^{\pm}$}
Let $(\pi, \mc P(i \lambda, +))$ be a unitarized principal series. Let $f \in \mc P(i \lambda, +)^{\infty}$ in the noncompact model. Fix the Whitakker functional $\exp -2 \pi i x \in \mc P(-i \lambda, +)^{-\infty}$. We have for $a>0$,
$$Wh_f(kan)=\langle \pi(k^{-1}) f, \pi^*(an_t) \exp(-2 \pi ix) \rangle=\exp 2 \pi i t \langle \pi(k^{-1}) f, a^{-1+ i \lambda}  \exp(-2 \pi i a^{-2} x) \rangle$$
$$=  a^{-1+ i \lambda} \exp 2 \pi i t \mc F(\pi(k)^{-1}f)(a^{-2}).$$
and 
 $$Wh_f(kI_{-}an)=\langle \pi( k^{-1}) f, \pi*(I_{-}) \pi^*(an_t) \exp(-2 \pi ix) \rangle$$
 $$=\exp 2 \pi i t \langle (\pi(k^{-1}) f)(x), a^{-1+ i \lambda} \pi^*(I_{-1}) \exp(-2 \pi i a^{-2} x) \rangle$$
 $$= \chi(-1) a^{-1+ i \lambda} \exp 2 \pi i t \mc F(\pi(k)^{-1}f)(-a^{-2}).$$
 To summarize, we have
 \begin{lem}\label{wh} In $SL^{\pm}(2)$, we have
 $$Wh_f(kan)=\chi(a) |a|^{-1+ i \lambda} \exp 2 \pi i t \mc F(\pi(k)^{-1}f)(a^{-2} \sgn(a)).$$
 \end{lem}

\subsection{Expression of $\int Wh_{f_1}(a) \overline{Wh_{f_2}(a)} \chi_0(a) |a|^{2-2s}  \frac{d a}{|a|}$}
Following Miller-Schmid \cite{ms0}, let 
$$\mc F(|x|^{-1+ u})(\xi)= G_0(u)|\xi|^{-u}, \qquad \mc F(\sgn(x) |x|^{-1+u})(\xi)=G_1(u) \sgn(\xi) |\xi|^{-u}$$ with
$$G_0(u)=2 (2\pi)^{-u} \Gamma(u) \cos(\frac{\pi u}{2}), \qquad G_1(u)=2 i (2 \pi)^{-u} \Gamma(u) \sin \frac{\pi u}{2}.$$

\noindent
Since the Fourier operator $\mc F$ is unitary and transforms multiplication into convolution, formally we have
\begin{equation}
\begin{split}
  & \int Wh_{f_1}(a) \overline{Wh_{f_2}(a)} \chi_0(a) |a|^{2-2s}  \frac{d a}{|a|}\\
  = & \int_{-\infty}^{\infty} Wh_{f_1}(a) \overline{Wh_{f_2}(a)}  \chi_0(a)|a|^{2-2s} \frac{d a}{|a|} \\
 = & \int_{-\infty}^{\infty} \chi_1(a) |a|^{-1+ i \lambda_1}  \mc F(f_1)(a^{-2} \sgn(a)) \overline{\chi_2(a)|a|^{-1+ i \lambda_2} \mc F(f_2)(a^{-2} \sgn(a))} \chi_0(a) |a|^{2-2s} \frac{ d a}{|a|} \\
 =& \frac{1}{2} \int_{-\infty}^{\infty} \chi_1(x) \chi_2(x) \chi_0(x) |x|^{-\frac{i\lambda_1-i \overline{\lambda_2}}{2}+s-1} 
 \mc F(f_1)(x) \overline{ \mc F(f_2)(x)} d x \\
 = & \frac{G_0(u)}{2} \int_{-\infty}^{\infty} \int_{-\infty}^{\infty} |x-y|^{-u} f_1(x) \overline{f_2(y)} d x dy \,\,\, \qquad (\chi_0\chi_1 \chi_2(-1)=1) \\
  &  \frac{G_1(u)}{2} \int_{-\infty}^{\infty} \int_{-\infty}^{\infty} \sgn(y-x) |y-x|^{-u} f_1(x) \overline{f_2(y)} d x dy \,\,\, \qquad (\chi_0\chi_1 \chi_2(-1)=-1)
\end{split}
\end{equation}
with $u=s-\frac{i\lambda_1-i \overline{\lambda_2}}{2}$.\\
\\
Due to the structure of this integral,  we can see that if $  \chi_1 \chi_2(-1)=1$, we have $L(s, \mc H_{\pi_1} \times \mc H_{\pi_2})=\sum_{n \neq 0} a_n \overline{b_n} |n|^{-s}$. If $\chi_1 \chi_2(-1)=-1$, we have  $L(s, \mc H_{\pi_1} \times \mc H_{\pi_2})=\sum_{n \neq 0} a_n \overline{b_n} \sgn(n) |n|^{-s}$.
If we use $E(s, \phi, g)_{o}$ to define the $L$-function, then it will be the other way around. Now we  choose the auxiliary actions of $\pi_1(I_{-})$ and $\pi_2(I_{-})$ to be even.  We have
\begin{lem}
Formally, for any $\phi \in C^{\infty}(K/M), f_i\in \mc P(i \lambda_i, +)^{\infty}_{e}, i \in [1,2]$, we have
$$Tr_{\mathbb R}^{RS}(f_1, f_2, (\phi_s)_{e})= \frac{G_0(u)}{2} \int \int_{-\infty}^{\infty} \int_{-\infty}^{\infty} |x-y|^{-u} (\pi_1(k)^{-1}f_1)(x) \overline{(\pi_2(k)^{-1}f_2)(y)} \phi(k)d x dy d k;$$
$$Tr_{\mathbb R}^{RS}(f_1, f_2, (\phi_s)_o)= \frac{G_1(u)}{2} \int \int_{-\infty}^{\infty} \int_{-\infty}^{\infty} \sgn(y-x) |y-x|^{-u} (\pi_1(k)^{-1}f_1)(x) \overline{(\pi_2(k)^{-1}f_2)(y)} \phi(k)d x dy d k.$$

\end{lem}

\subsection{Expression of $\Delta_{+}(f_1, f_2, u)$}
Write 
$$\Delta_{+}(f_1, f_2, u)=\int_{-\infty}^{\infty} \int_{-\infty}^{\infty} |x-y|^{-u} f_1(x) \overline{f_2(y)} d x dy.$$
We will convert this integral into an integral in the compact picture. 
We have
$$f_1(x)=(1+x^2)^{-\frac{1+ i \lambda_1}{2}}\tilde{f_1}(\cot^{-1} x), \qquad \theta_1=\cot^{-1} x;$$
$$f_2(y)=(1+y^2)^{-\frac{1+ i \lambda_2}{2}}\tilde{f_2}(\cot^{-1} y), \qquad \theta_2=\cot^{-1} y$$
$$\tilde f_i(\pi+\theta)=\tilde f_i(\theta).$$
Hence
$$\Delta_{+}(f_1, f_2, u)=\int_{0}^{\pi} \int_{0}^{\pi}
|\cot \theta_1-\cot {\theta_2}|^{-u} \tilde{f_1}(\theta_1) \overline{\tilde f_2}(\theta_2) |\sin \theta_1|^{-1+ i \lambda_1}|\sin \theta_2|^{-1- i \overline{\lambda_2}} d \theta_1 d \theta_2 $$
$$=\frac{1}{4}\int_{0}^{2 \pi} \int_{0}^{2\pi}
|\sin( \theta_1-\theta_2)|^{-u} \tilde{f_1}(\theta_1) \overline{\tilde f_2}(\theta_2) |\sin \theta_1|^{-1+u+ i \lambda_1}|\sin \theta_2|^{-1+u- i \overline{\lambda_2}} d \theta_1 d \theta_2.$$
\begin{thm}\label{rsr}
Let $\lambda_i \in \mathbb R$. Then formally
$$Tr_{\mathbb R}^{RS}(f_1, f_2, \phi_s)=\frac{G_0(u)}{8}\int_{0}^{2 \pi} \int_{0}^{2\pi} |\sin( \theta_1-\theta_2)|^{-u}  |\sin \theta_1|^{-1+u+ i \lambda_1}|\sin \theta_2|^{-1+u- i \lambda_2} H_{f_1, f_2,\phi}(\theta_1, \theta_2) d \theta_1 d \theta_2,$$
where $u=s-\frac{i\lambda_1-i \lambda_2}{2}$ and
$$H_{f_1, f_2, \phi}(\theta_1, \theta_2)=\frac{1}{2\pi} \int \phi(\theta) \tilde{f_1}(\theta_1+\theta) \overline{\tilde f_2}(\theta_2+\theta) d \theta.$$
If $\Re(s) \in (0,1)$, these integrals converge absolutely for all $f_i \in \mc P(i \lambda_i, +)$ and $\phi \in C^{\infty}(K/M)$. Hence, we obtain a continuous $G$-equivariant trilinear form
$$Tr_{\mathbb R}^{RS}: \mc P(i \lambda_1,+)_{e} \times \mc P(i \lambda_2, +)_{e} \times \mc P(2s-1,+)^{\infty}_{e} \rightarrow \mathbb C.$$
\end{thm}
Proof:  The integral
$$\int_{K} \pi_1(k^{-1}) \tilde f_1(\theta_1) \overline{\pi_2(k^{-1}) \tilde f_2(\theta_2)} \phi_(k) d k= \frac{1}{2\pi} \int \phi(\theta) \tilde{f_1}(\theta_1+\theta) \overline{\tilde f_2}(\theta_2+\theta) d \theta.$$
We obtain the desired expression of $Tr_{\mathbb R}^{RS}(f_1, f_2, \phi_s)$. For unitary principal series, $\tilde f_i \in L^2(K/M)$. Since $\phi$ is smooth, $H_{f_1,f_2,\phi}(\theta_1, \theta_2)$ is continuous and bounded. Our theorem follows from the absolute convergence
of $$\int_{0}^{2 \pi} \int_{0}^{2\pi}
 |\sin( \theta_1-\theta_2)|^{-\Re(u)}  |\sin \theta_1|^{-1+\Re(u)}|\sin \theta_2|^{-1+\Re(u)} d \theta_1 d \theta_2,$$
 when $\Re(u) \in (0,1)$.
 $\Box$\\
 \\
 From our proof, there is enough space to perturb the parameter $ \lambda_1, \lambda_2$ to make the integral converges even for the complementary series. However, for complementary series, the Hilbert norm $\| \tilde f \|$ is no longer the $L^2$-norm on $K$. The computation will be complicated and less illuminating. In addition, by the Selberg's conjecture, complementary series will not occur in $L^2(G/\Gamma)$ for a congruence subgroup $\Gamma$.  We will exclude the complementary series in the next two sections. \\
 \\
Finally, $\Delta_{-}(f_1, f_2, \phi_s)$ can be defined and computed similarly. Similar theorem can be proved for
$$Tr_{\mathbb R}^{RS}: \mc P(i \lambda_1,+)_{e} \times \mc P(i \lambda_2, +)_{e} \times \mc P(2s-1,+)^{\infty}_{o} \rightarrow \mathbb C.$$
We leave these computation and verification to the reader. The L-function corresponding to this computation is $\sum_{n \neq 0} \sgn(n) a_n \overline{b_n} |n|^{-s}$ with $E(s, *,*)_{o}$. In the next section, we shall stick with $E(s,*,*)_{e}$ and 
$L(s, \mc H_{\pi_1} \times \mc H_{\pi_2})=\sum_{n \neq 0}  a_n \overline{b_n} |n|^{-s}$.
\section{ Bounds for $\mc I_{\pi_1, \pi_2}(s)$: Spherical Unitary Principal series of $SL^{\pm}(2)$} 
Let $(\pi_1, \pi_2)=(\mc P(i \lambda_1, +)_{e}, \mc P(i \lambda_2, +)_{e})$.
By Cor. \ref{main3}, in order to bound $L(s, \mc H_{\pi_1} \times \mc H_{\pi_2})=\sum_{n \neq 0}  a_n \overline{b_n} |n|^{-s}$ in the $s$-aspect, we must find lower bound for  $\mc I_{\pi_1, \pi_2}(s)$. The index $\mc I_{\pi_1, \pi_2}(s)$ is defined to be
$$\sup \{\frac{|Tr_{\mathbb R}^{RS}(f_1, f_2, {\mathbf 1}_s)|}{\|f_1\| \|f_2\|} \,\,  : \,\,\, \|f_i \| \neq 0 \,\,\, i=0,1 \}.$$
Since
$$H_{f_1,f_2 \mathbf 1}(\theta_1, \theta_2)=\frac{1}{2\pi} \int  \tilde{f_1}(\theta_1+\theta) \overline{\tilde f_2}(\theta_2+\theta) d \theta,$$ 
this function is $K$-invariant. We may write
$$ H_{f_1,f_2, \mathbf 1}(\theta_1, \theta_2)=H^0_{f_1,f_2}(\theta_1-\theta_2).$$
Here we assume that $\theta_1, \theta, \theta_2 \in \mathbb R/ 2 \pi \mathbb Z$.
Then for $Re(s) \in (0,1)$, by Theorem \ref{rsr}, we have
\begin{equation}
\begin{split}
 & Tr_{\mathbb R}^{RS}(f_1, f_2, \mathbf 1_{s}) \\
=& \frac{G_0(u)}{8}\int_{0}^{2 \pi} \int_{0}^{2\pi} |\sin( \theta_1-\theta_2)|^{-u}  |\sin \theta_1|^{-1+u+ i \lambda_1}|\sin \theta_2|^{-1+u- i \lambda_2} H_{f_1, f_2}^0(\theta_1- \theta_2) d \theta_1 d \theta_2\\
= & \frac{G_0(u)}{8}\int_{0}^{2 \pi} \int_{0}^{2\pi} |\sin( \theta_1-\theta_2)|^{-u}  |\cos \theta_1|^{-1+u+ i \lambda_1}|\cos \theta_2|^{-1+u- i \lambda_2} H_{f_1, f_2}^0(\theta_1- \theta_2) d \theta_1 d \theta_2 \\
= &  \frac{G_0(u)}{8}\int_{0}^{2 \pi} \int_{0}^{2\pi} |\sin( \theta_1)|^{-u}  |\cos (\theta_1+\theta_2)|^{-1+u+ i \lambda_1}|\cos \theta_2|^{-1+u- i \lambda_2} H_{f_1, f_2}^0(\theta_1) d \theta_1 d \theta_2\\
= & \frac{G_0(u)}{8}\int_{0}^{2 \pi} \int_{0}^{2\pi} |\sin( \theta_1)|^{-u}  K(\theta_1) H_{f_1, f_2}^0(\theta_1) d \theta_1 d \theta_2.
\end{split}
\end{equation}
Here $K(\theta_1)=\int |\cos (\theta_1+\theta_2)|^{-1+u+ i \lambda_1}|\cos \theta_2|^{-1+u- i \lambda_2} d \theta_2$.\\
\\
Since $H$ and $K$ all appear to be  convolutions, we shall apply the theory of Fourier series. 
The Fourier coefficients will be denoted by
$$\hat{\tilde f}(n)=\frac{1}{2\pi} \int_{0}^{2 \pi} \tilde f(\theta) \exp( - i n \theta )d \theta.$$
 Then for any $\tilde f \in C^{\infty}(K)$, we have
$$\tilde f(\theta)=\sum_{n} \hat{\tilde f}(n) \exp (i n \theta).$$
We use $\| \hat{\tilde f} \|_{L^1}$ to denote the $L^1$-norm of the Fourier coefficients of $\tilde f$,  $\| \hat{\tilde f} \|_{L^2}$ to denote the $L^2$-norm, and $\| \hat{\tilde f} \|_{\sup}$ to denote the sup-norm.
\begin{lem}\label{l1l2} For $(\pi_1, \pi_2)=(\mc P(i \lambda_1,+)_{e}, \mc P(i \lambda_2, +)_{e})$ with $\lambda_1, \lambda_2 \in \mathbb R$, we have
$$I_{\pi_1, \pi_2}(s) = \sup \{\frac{|G_0(u) \int H(\theta) K(\theta) |\sin \theta|^{-u} d \theta|}{8 \pi \|\hat H\|_{L^1}}: \|\hat H\|_{L^1} < \infty \}$$
with $K(\theta)=\int |\cos (\theta+\theta_2)|^{-1+u+ i \lambda_1}|\cos \theta_2|^{-1+u- i \lambda_2} d \theta_2$ and $u=s-\frac{i\lambda_1-i \lambda_2}{2}$.
\end{lem}
Proof: Notice that 
$$H_{f_1, f_2}^0(\theta)= \frac{1}{2 \pi} \int \tilde{f_1}(\theta+\theta_1) \overline{\tilde{f_2}} (\theta_1) d \theta_1$$
$$= \sum \widehat{\tilde{ f_1}}(n) \overline{\widehat{\tilde f_2}}(n) \exp i n \theta$$
By Cauchy-Schwartz inequality,
$$\| f_1\| \|f _2 \| = \pi \| \tilde f_1 \| \|\tilde f_2 \|=\pi \|\widehat{\tilde f_2} \|_{L^2} \| \widehat{\tilde f_2} \|_{L^2}\geq \pi \| \widehat{H^0_{f_1, f_2}} \|_{L^1} $$
The equality can hold for proper choices of $\tilde f_1$ and $\tilde f_2$. Hence
$$\frac{Tr_{\mathbb R}^{RS}(f_1, f_2, \mathbf 1)}{\|f_1\| \| f_2\|} \leq \frac{|G_0(u) \int H_{f_1, f_2}^0 (\theta) K(\theta) |\sin \theta|^{-u} d \theta|}{8 \pi \|\widehat{H^0_{f_1,f_2}}\|_{L^1}}.$$
with equality for a proper choice of $(f_1, f_2)$. Our lemma follows immediately.
$\Box$ 
\begin{thm}\label{isup}
For $(\pi_1, \pi_2)=(\mc P(i \lambda_1,+)_{e}, \mc P(i \lambda_2, +)_{e})$ with $\lambda_1, \lambda_2 \in \mathbb R$,  we have
$$\mc I_{\pi_1, \pi_2}(s) =  \frac{|G_0(u)| \| \widehat{ K(\theta) |\sin \theta|^{-u}} \|_{\sup}}{8 \pi }$$
with $K(\theta)=\int |\cos (\theta+\theta_2)|^{-1+u+ i \lambda_1}|\cos \theta_2|^{-1+u- i \lambda_2} d \theta_2$ and $u=s-\frac{i\lambda_1-i \lambda_2}{2}$.
\end{thm}
Proof: By Lemma \ref{l1l2}, we will need to find
$$\sup \{ \frac{|\int H(\theta) K(\theta) |\sin \theta|^{-u} d \theta|}{\|\hat{H}(n)\|_{L^1}} \}=\sup \{\frac{|\sum_{n \in \mathbb Z} \hat H(n)  \widehat{ K(\theta) |\sin \theta|^{-u}}(-n)| }{\| \hat{H}(n) \|_{L^1}}  \}.$$
Clearly 
$$\sup \{\frac{|\sum_{n \in \mathbb Z} \hat H(n)  \widehat{ K(\theta) |\sin \theta|^{-u}}(-n)| }{\| \hat{H}(n) \|_{L^1}}  \}= \sup \{ |\widehat{ K(\theta) |\sin \theta|^{-u}}(n)|: n \in \mathbb Z \}.$$
Our theorem then follows. $\Box$\\
\\
The Fourier coefficients of $K(\theta) |\sin \theta|^{-u}$ can be expressed explicitly in an infinite sum, as a convolution of two sequences. But they are difficult to compute. At least, we can prove the following
\begin{thm}\label{estimatei}Let $s=s_0+ i s_1$.
For unitary principal series $(\pi_i, \mc P(i \lambda_i, +)_{e})$, there exists constant $C_{\lambda_1, \lambda_2, s_0}$ such that
$$I_{\pi_1, \pi_2}(s) \geq C_{\lambda_1, \lambda_2, s_0} (1+ |s_1|)^{s_0-1}$$ for every $1> s_0 >\frac{1}{2}$.
\end{thm}
There is a natural barrier at $s_0=\frac{1}{2}$ in our proof. It is likely that this bound can be improved, perhaps to the Weyl bound, $(1+|s_1|)^{-\frac{1}{3}+\epsilon}$ for $s_0$ sufficiently close to $\frac{1}{2}$. \\
\\
Proof:  Since $\|\hat{H} \|_{L^1} < \infty$, $H(\theta)$ must be bounded and continuous. Let $\epsilon > 0 $. We  choose the test function $H(\theta)=\overline{K(\theta)} |\sin \theta|^{\epsilon+ i \Im(u)}$ to make $H(\theta) K(\theta) |\sin \theta |^{-u}$ nonegative  and not vary with $s_1$. Hence 
$$\int H(\theta) K(\theta) |\sin \theta |^{-u} d \theta= C_{s_0, \epsilon}^{(1)} >0.$$
is fixed for all $s_1$. To apply Lemma \ref{l1l2}, we must estimate $\|\hat{H} \|_{L^1}$. \\
\\
Since $K(\theta)=\int | \cos(\theta+\theta_2)|^{-1+s+\frac{i \lambda_1+i \lambda_2}{2}} | \cos(\theta_2)|^{-1+s-\frac{i \lambda_1+ i \lambda_2}{2}} d \theta_2$, 
we have
$$\hat K(n)= \widehat{|\cos(\theta)|^{-1+s-\frac{i \lambda_1+ i \lambda_2}{2}}} (n)  \widehat{ | \cos(\theta)|^{-1+s+\frac{i \lambda_1+i \lambda_2}{2}}} (n) .$$
By Cauchy-Schwartz inequality, 
$$\|\widehat{\overline{K(\theta)}} \|_{L^1}=\|\widehat{K(\theta)} \|_{L^1} \leq \| \widehat{|\cos(\theta)|^{-1+s-\frac{i \lambda_1+ i \lambda_2}{2}}} \|_{L^2} \| \widehat{ | \cos(\theta)|^{-1+s+\frac{i \lambda_1+i \lambda_2}{2}}} \|_{L^2} $$
$$= \||\cos(\theta)|^{-1+s-\frac{i \lambda_1+ i \lambda_2}{2}} \|_{L^2} \|  | \cos(\theta)|^{-1+s+\frac{i \lambda_1+i \lambda_2}{2}} \|_{L^2}= \| |\cos(\theta)|^{-1+s_0} \|_{L^2}^2=C_{s_0}^{(2)}.$$
Since $H(\theta)=\overline{K(\theta)}|\sin \theta|^{\epsilon+ i s_1-\frac{ i \lambda_1-i \lambda_2}{2}}$,
Its Fourier coefficients $\hat H$ will be a convolution of $\widehat{\overline{K}}$ and $\widehat{|\sin \theta|^{\epsilon+ i s_1-\frac{ i \lambda_1-i \lambda_2}{2}}}$. We have
$$\| \hat{H} \|_{L^1} \leq \|\hat{K}\|_{L^1} \| \widehat{|\sin \theta|^{\epsilon+ i s_1-\frac{ i \lambda_1-i \lambda_2}{2}}}\|_{L^1} \leq C_{s_0}^{(2)} \| \widehat{|\sin \theta|^{\epsilon+ i s_1-\frac{ i \lambda_1-i \lambda_2}{2}}}\|_{L^1}. $$
We will show in Lemma \ref{l1} that  $\| \widehat{|\sin \theta|^{\epsilon+ iu_1}}\|_{L^1} \leq C_{\epsilon} (1+|u_1|)^{\frac{1}{2}}$ with $u_1= s_1-\frac{  \lambda_1- \lambda_2}{2}$. 
Our theorem then follows from Lemma \ref{l1l2} and $|G_0(s_0+i s_1)| \geq C^{(3)}_{s_0} (1+ |s_1|)^{s_0-\frac{1}{2}}$. $\Box$
\begin{lem}\label{l1} For any $\epsilon>0$ and $u_1 \in \mathbb R$, we have
$$\| \widehat{|\sin \theta|^{\epsilon+ iu_1}}\|_{L^1} \leq C_{\epsilon} (1+|u_1|)^{\frac{1}{2}}.$$
\end{lem}
Proof: Let $p(\theta)=|\cos (\theta)|^{\epsilon+ iu_1}$. Observe that $\widehat{|\sin \theta|^{\epsilon+ iu_1}}(m)$ differs from $\hat{p}(m)$ by a phase factor $\exp \frac{i \pi m}{2}$. It suffices to prove that $\| \hat{p} \|_{L^1} \leq C_{\epsilon} (1+|u_1|)^{\frac{1}{2}}.$
To compute the  Fourier coefficients of $p(\theta)$, we apply the following well-known formula
(\cite{pk})
$$\int_{-\frac{\pi}{2}}^{\frac{\pi}{2}} (\cos t)^{a+b-2} \exp i(a-b+2x) t d t= \frac{\pi \Gamma(a+b-1)}{2^{a+b-2} \Gamma(a+x) \Gamma(b-x)}.$$ 
Then
$$\hat{p}(2m)=\frac{1}{2 \pi} \int_{0}^{2\pi} p(\theta) \exp 2mi \theta d \theta= \frac{\Gamma( \epsilon+ i u_1+1)}{2^{\epsilon+ i u_1} \Gamma(m+ \frac{\epsilon+ i u_1}{2}+1) \Gamma(-m+1+ \frac{\epsilon+ i u_1}{2})}.$$
Obviously $\hat{p}(2m)=\hat{p}(-2m)$ and $\hat{p}(2m+1)=0$. Applying  the duplication formula of $\Gamma$ function, we obtain
$$\hat{p}(2m)=\frac{\Gamma(\epsilon+ i u_1+1) \Gamma(m- \frac{\epsilon+ i u_1}{2}) (-1)^{m+1} \sin( \frac{\pi \epsilon}{2}+ \frac{i \pi u_1}{2})}{2^{\epsilon+ i u_1} \pi \Gamma(m+\frac{\epsilon+ i u_1}{2}+1)}.$$
The factor 
$$|\frac{\Gamma(\epsilon+ i u_1+1)  \sin( \frac{\pi \epsilon}{2}+ \frac{i \pi u_1}{2})}{2^{\epsilon+ i u_1} \pi }| \leq C_{\epsilon}^{(1)} (1+ | u_1|)^{\epsilon+ \frac{1}{2}}.$$
In addition,
$$|\frac{\Gamma(-\frac{\epsilon+ i u_1}{2})}{\Gamma(\frac{\epsilon+ i u_1}{2}+1)}|+ 2 \sum_{m=1}^{\infty} |\frac{\Gamma(m-\frac{\epsilon+ i u_1}{2})}{\Gamma(m+\frac{\epsilon+ i u_1}{2}+1)}| \leq C_{\epsilon}^{(2)} (1+ |u_1|)^{-\epsilon}.$$
Our Lemma is then proved. We may simply fix $\epsilon=1$ when we apply this lemma. $\Box$ \\
\\
Perhaps, it is not difficult to see how our computation of $\mc I_{\pi_1, \pi_2}(s)$ can be modified for the nonspherical unitary principal series $(\pi_1, \pi_2)$. One will have to insert certain $\sgn$-functions into the expressions in Theorem \ref{rsr}. One would also have to modify Theorem \ref{isup}. Theorem \ref{estimatei} is expected to be true for nonspherical unitary principal series.

\commentout{\section{$Tr_{\mathbb R}^{RS}$: nonspherical principal series case}
Let $u \in (-1, 0)$ and $\lambda \in \mathbb R$. The $K$-types in $\mc P(i \lambda, -)$ are
$$v_{2m+1}^{(i \lambda)}(x)=(1+x^2)^{\frac{-i \lambda+u}{2}} (\frac{1+xi}{1-xi})^{m+\frac{1}{2}} \qquad (m \in \mathbb Z).$$
Here $x= \cot \theta \,\,\,(\theta \in (0, \pi))$ and $(\frac{1+xi}{1-xi})^{m+\frac{1}{2}}= \exp i (2m+1) \theta$ is well-defined.}
\section{ Growth of $L$-function}
By Theorems \ref{estimatei} and \ref{bound1}, we have 
\begin{thm}
Let $\Gamma$ be a nonuniform lattice of $SL(2)$ such that $\Gamma \cap N=N_{\mathbb Z}$. Let 
$(\pi_1, \pi_2)=(\mc P(i \lambda_1, +)_{e}, \mc P(i \lambda_2, +)_{e})$. Let $\mc H_{\pi_1}$ and 
$\mc H_{\pi_2}$ be two cuspidal automorphic representation of $SL^{\pm}(2)/\Gamma$. We have
for  $1> s_0 > \frac{1}{2}$
$$L(s_0+i s_1, \mc H_{\pi_1} \times \mc H_{\pi_2}) \leq  C_{v, s_0, \mc H_{\pi_1}, \mc H_{\pi_2}} \sup_G(|E(s, \phi_s, g) v(g)|)
(1+\|s_1\|)^{1-s_0}  $$
where $v(g)$ is a bounded continuous positive function on $\mc F$ such that $v(ka_i n_i)=|a_i|^{2-\epsilon}$ for a fixed  $\epsilon \in (0, \min(2\Re(s), 2- 2\Re(s))$ on the every Siegel set $S_i$.
\end{thm}
Applying the bounds on the Eisenstein series due to Young, Huang-Xu, Assing and Nordentoft, we have
\begin{thm}
Let $\Gamma$ be a congruence subgroup of $SL(2, \mathbb Z)$ or a     conjugate of such a group with $\Gamma \cap N=N_{\mathbb Z}$. Let 
$(\pi_1, \pi_2)=
(\mc P(i \lambda_1, +)_{e}, \mc P(i \lambda_2, +)_{e})$. Let $\mc H_{\pi_1}$ and 
$\mc H_{\pi_2}$ be two cuspidal automorphic representation of $SL^{\pm}(2)/\Gamma$. Then for $s_0 > \frac{1}{2}$ and $\epsilon>0$
$$L(s_0+i s_1, \mc H_{\pi_1} \times \mc H_{\pi_2}) \leq  C_{ s_0, \epsilon, \mc H_{\pi_1}, \mc H_{\pi_2}} 
(1+\|s_1\|)^{\frac{11}{8}-s_0+\epsilon}  $$
If $\Gamma=SL(2, \mathbb Z)$, then
$$L(s_0+i s_1, \mc H_{\pi_1} \times \mc H_{\pi_2}) \leq  C_{ s_0, \epsilon, \mc H_{\pi_1}, \mc H_{\pi_2}} 
(1+\|s_1\|)^{\frac{4}{3}-s_0+\epsilon}  $$ 
\end{thm} 
For $SL(2, \mathbb Z)$, $L$-function satisfies a functional equation. By  Phragmen-Lindelof principle, we have for any $t \in \mathbb R$,
$$L(\frac{1}{2}+ i t, \mc H_{\pi_1} \times \mc H_{\pi_2}) \leq  C_{ \epsilon, \mc H_{\pi_1}, \mc H_{\pi_2}} 
(1+\|t\|)^{\frac{5}{6}+\epsilon}.  $$
This breaks the convexity bound. Recall that  $L(\frac{1}{2}+ t, \mc H_{\pi_1} \times \mc H_{\pi_2})$ is a degree 4 $L$-function. Its convexity bound is
$$L(\frac{1}{2}+ i t, \mc H_{\pi_1} \times \mc H_{\pi_2}) \leq  C_{ \epsilon, \mc H_{\pi_1}, \mc H_{\pi_2}} 
(1+\|t\|)^{1+\epsilon}.  $$
Our method does apply in the general setting of 
Rankin-Selberg $L$-function for $GL(n) \times GL(n)$ (\cite{go} \cite{cogdell}). 
In the adelic setting, $L(s, \mc \pi_1 \times \pi_2)$ will be a degree $n^2$ $L$-function and the convexity bound Of $L(\frac{1}{2}+ it)$ will be $(1+|t|)^{\frac{n^2}{4}+\epsilon}$ (\cite{sarnak}). Theorem \ref{main2} and Theorem \ref{main3} should hold for $GL(n)$. For $GL(n)$ the contribution from $\sup_{G}(|E(s, g) v(g)|)$ will be very small. Hence the main contribution to our bound 
$$
L(s, \mc H_{\pi_1} \times \mc H_{\pi_2}) | \leq C_{v, \mc H_{\pi_1}, \mc H_{\pi_2}} \frac{\sup_G(|E(s, g) v(g)|)}{\sup\{|Tr_{\mathbb R}^{RS}(f_1, f_2, {\mathbf 1}_s)|\,\,  : \,\,\, \|f_1\|=\|f_2\|=1 \} }$$
will be the Archimedean index 
$$\mc I_{\pi_1, \pi_2}(s)= \sup\{|Tr_{\mathbb R}^{RS}(f_1, f_2, {\mathbf 1}_s)|\,\,  : \,\,\, \|f_1\|=\|f_2\|=1 \}.$$
Precise value of this index will be a challenging problem in representation theory. Nevertheless,
we expect $\mc I_{\pi_1, \pi_2}(s)$ can be similarly estimated using the method in this paper. It is likely that one will be able to establish the subconvexity bound for all Rankin-Selberg function $L(\frac{1}{2}+ it)$ for $GL(n) \times GL(n)$ in a similar fashion.

\end{document}